\documentclass[10pt,english]{article}
\usepackage{mathptmx}
\usepackage[T1]{fontenc}
\usepackage[latin9]{luainputenc}
\usepackage{babel}
\usepackage{amsmath}
\usepackage{amssymb}
\usepackage{graphicx}
\usepackage{setspace}
\usepackage[pdfusetitle,
 bookmarks=true,bookmarksnumbered=false,bookmarksopen=false,
 breaklinks=false,pdfborder={0 0 1},backref=false,colorlinks=false]
 {hyperref}
\begin{document}
\title{Elephant Random Walk with multiple extractions}
\author{Simone Franchini{\normalsize\thanks{Correspondence: simone.franchini@yahoo.it}\thanks{Sapienza Università di Roma, Piazza Aldo Moro 1, 00185 Roma, Italy}}}
\date{~}
\maketitle
\begin{abstract}
Consider a generalized Elephant Random Walk in which the step is chosen
by selecting $k$ previous steps with $k$ odd and then going in the
majority direction with a probability $p$ and in the opposite direction
otherwise. In the $k=1$ case the model is the original one and could
be resolved exactly by analogy with Friedman's urn. However the analogy
cannot be extended to the $k>2$ case already. In this paper we show
how to treat the model for each $k$ by analogy with the more general
urn model of Hill, Lane and Sudderth. Interestingly for $k>2$ we
found a critical dependence from the initial conditions beyond a certain
values of the memory parameter $p$, and regions of convergence with
entropy that is sub-linear in the number of steps.

~

\noindent\textit{keywords: urn models, increasing returns, stochastic
approximation, lattice field theory}

~

~

~

~

~

~

~

~

~

~

~

~

~

~
\end{abstract}
\pagebreak{}

\section{Introduction\protect\label{sec:Introduction}}

The Elephant Random Walk (ERW) is a simplified model of a random walk
with long-range memory, in which each step is determined by extracting
one of the previous, then going in the same direction with probability
$p$, and in the opposite otherwise. \cite{ERW=000020shcutz=000020trimper,Gut_Stadmuller,Bercu,Maulik}
This simple model has received much attention since its introduction
and is still a very active area, see e.g. \cite{Nakano,Hu-Fan,Akimoto,Dhillon,Atreiu,Quin}.

An important advancement in the study of this model has been made
in 2016, when Baur and Bertoin noticed in \cite{ERW=000020UM=000020Baur=000020Berton}
that the ERW can be mapped exactly into a two colors urn model of
the Friedman type \cite{Pemantle=000020review,Mahmoud}, where at
each step a ball is drawn from the urn and then replaced, together
with a finite number of new balls, possibly of both colors, depending
on which of ball was drawn. 

This finding allowed many quantities of interest to be computed from
known results on these types of urns. Unluckily the analogy with the
Friedman's type urns cannot be extended to more complex memory mechanisms
where the step is chosen at the best of multiple draws. Let us consider
a generalized version of the ERW, in which the step is chosen at the
best of $k$ draws \cite{Dosi=000020Ermoliev}: define 
\begin{equation}
X:=\{X_{1},\,X_{2},\,...\,,X_{N}\}
\end{equation}
with $X_{n}\in\left\{ +,-\right\} $ for all $1\leq n\leq N$. The
average step up to time $N$ is defined as 
\begin{equation}
x_{N}:=\frac{1}{N}\sum_{n\leq N}X_{n}.
\end{equation}
Then, the next step $X_{N+1}$ is determined as follows: first, extract
$k$ previous steps from $X$, with $k$ odd integer. If the sum of
these steps is positive (means there is a majority of positive steps)
then $X_{N+1}$ is positive with probability $p$ and negative otherwise.
In case the sum is negative, then $X_{N+1}$ is negative with probability
$p$, and positive otherwise. For $k=1$ the original ERW is recovered.
In the present work we show how the analogy with a more general non-linear
urn, known as the Hill, Lane and Sudderth (HLS) urn model, \cite{HLS1,HLS2,Pemantle=000020Touch,Franchini,FB,FranchiniHLS2025}
(for which the Friedman urn represents a linear subcase and a mathematical
treatment for the small and large deviations exists already \cite{Franchini,FB,FranchiniHLS2025})
allows to deal with such generalizations too. Interestingly, for multiple
extractions there is a value $p_{c}$ of the memory parameter above
which the limit
\begin{equation}
x:=\lim_{N\rightarrow\infty}x_{N},
\end{equation}
critically depends on the initial condition $x_{M}$, fixed at some
finite time $M<\infty$. We also study the limit entropy density for
$x_{N}$ converging to some given $x$, 
\begin{equation}
\rho\left(x\right):=\lim_{N\rightarrow\infty}\frac{1}{N}\log P\left(x_{N}=\left\lfloor xN\right\rfloor /N\right),
\end{equation}
finding that for $k>2$ and above $p_{c}$ there is a region in the
$p$ vs $x$ plane where $\rho\left(x\right)=0$, i.e. where the entropy
is sub linear in the number of steps. 
\begin{figure}
\begin{centering}
\includegraphics[scale=0.29]{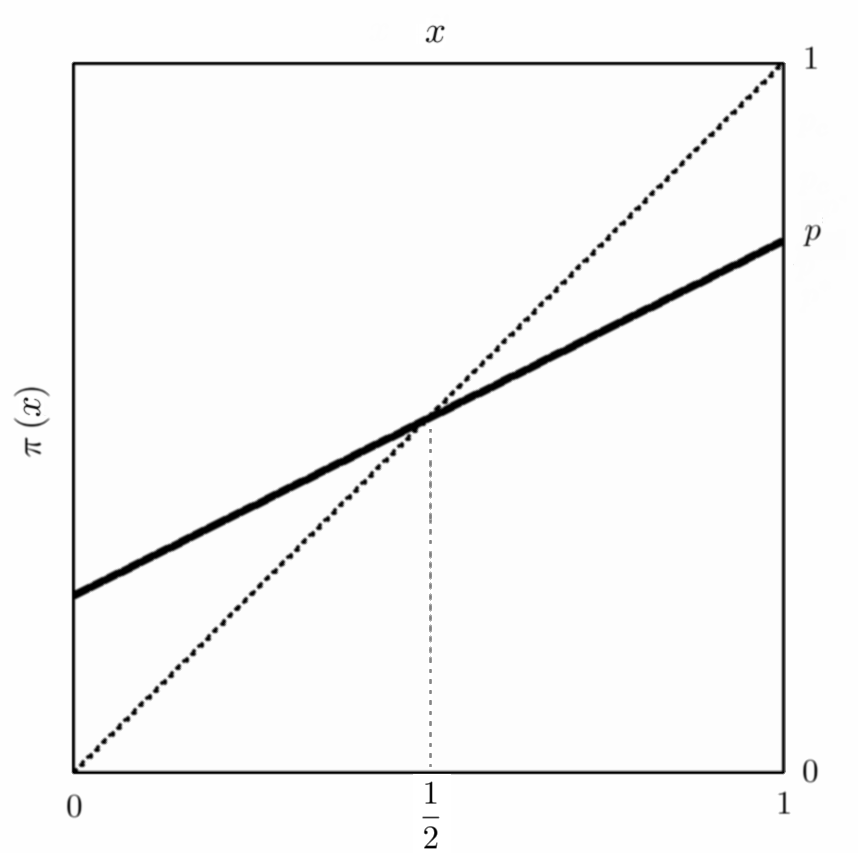}
\par\end{centering}
\caption{\protect\label{fig1}Example of linear urn function $\pi_{1}\left(y\right)$
for the classic ERW with $k=1$, the memory parameter is $p^{*}=3/4$.
The urn function always crosses the diagonal at $y_{0}=1/2$. \cite{FB}}
\end{figure}
\begin{figure}
\centering{}\includegraphics[scale=0.29]{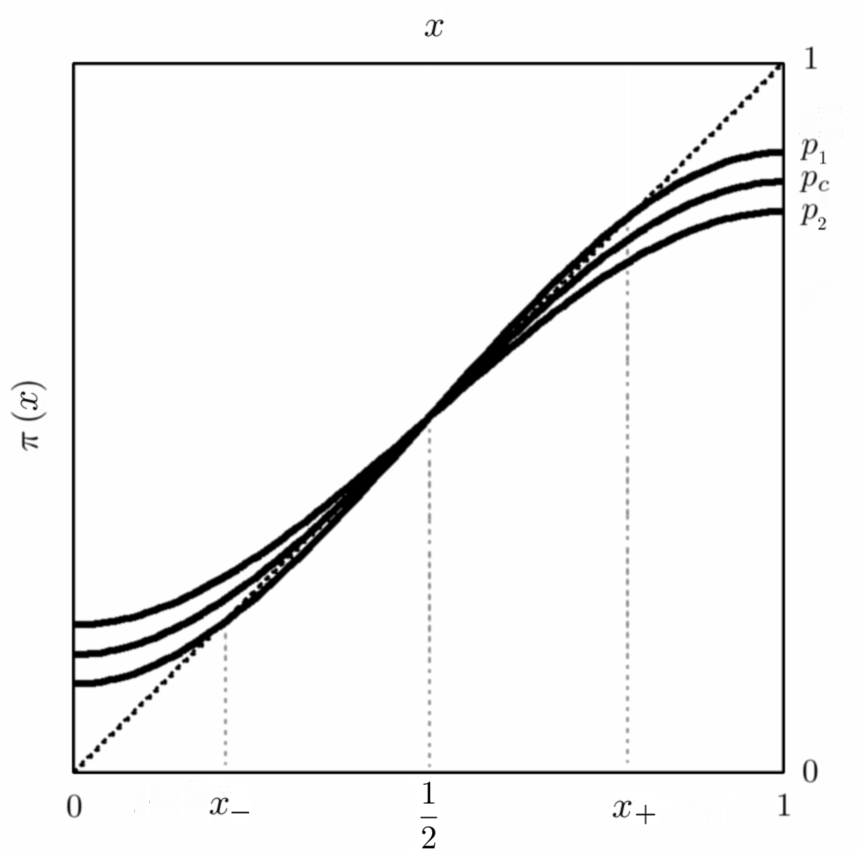}\caption{\protect\label{fig:2}Three examples of the the urn function $\pi_{3}\left(y\right)$
for the generalized ERW with $k=3$. The figure shows urn functions
for three non-trivial memory parameters, $p^{*}=2/3$, $p_{c}=5/6$
and $p^{**}=11/12$. Below $p_{c}$ the urn function down-crosses
the line $y$ at $y_{0}=1/2$, that is the only convergence point.
For $p>p_{c}$ the point $y_{0}$ becomes an up-crossing (unstable
equilibrium), and the urn function crosses the diagonal $y$ also
in $y_{-}$ and $y_{+}$, that are both down\textendash crossings
and are the new stable attractors for the process $y_{N}$. \cite{FB}}
\end{figure}

\section{Relation with HLS urns\protect\label{sec:Relation-with-HLS}}

An HLS urn model \cite{Pemantle=000020review,Mahmoud,HLS1,HLS2,Pemantle=000020Touch,Franchini,FB}
is a two color urn process governed by a functional parameter $\pi\left(y\right)$
called the \textit{urn function}. Let us consider an infinite capacity
two-color urn containing black and white balls, and let 
\begin{equation}
Y=\left\{ \,Y_{1},Y_{2},\,...,Y_{N}\right\} 
\end{equation}
be the process describing the type of ball $Y_{n}\in\left\{ 0,1\right\} $
that has been extracted at each time $n\leq N$. We associate the
value $1$ to black balls and $0$ to white balls. The process $Y$
evolves as follows: let the quantity
\begin{equation}
y_{N}:=\frac{1}{N}\sum_{n\leq N}Y_{t}
\end{equation}
be the density of black balls at time $N$, then at step $N+1$, a
new ball is added, black with probability $\pi\left(y_{N}\right)$
and white otherwise. The limit density is 
\begin{equation}
y:=\lim_{N\rightarrow\infty}y_{N}.
\end{equation}
The analogy between the ERW before and the HLS urns is easily unveiled
by taking $X_{n}=2Y_{n}-1$ and then writing the density of black
balls in terms of the average step 
\begin{equation}
y=\frac{1+x}{2}.\label{eq:trans}
\end{equation}

Start with $k=1$, ie the simple ERW: the probability of extracting
a black ball from the urn is $y$, then the ERW will go in the positive
direction with probability 
\begin{equation}
\pi_{1}\left(y\right):=p\,y+\left(1-p\right)\left(1-y\right)=\left(1-p\right)+\left(2p-1\right)y,\label{fdg}
\end{equation}
that is a linear urn function. \cite{Jack=000020Harris} In fact,
it is proven in \cite{Franchini} that linear HLS urns of the kind
$\pi\left(y\right)=a+by$ include the Friedman urn model as special
case. Let now move to the case $k=3$: the probability of taking a
positive step is that of extracting two positive and one negative,
plus that of extracting three positive: 
\begin{equation}
P_{3}\left(y\right):=y^{3}+3y^{2}\left(1-y\right)=3y^{2}-2y^{3},
\end{equation}
then, the corresponding urn function is
\begin{multline}
\pi_{3}\left(y\right):=p\,\left(3y^{2}-2y^{3}\right)+\left(1-p\right)\left(1-\left(3y^{2}-2y^{3}\right)\right)=\\
=\left(1-p\right)+\left(2p-1\right)\left(3y^{2}-2y^{3}\right)\label{eq:pik-1}
\end{multline}
and cannot be reduced to the linear case. In general, the probability
of finding a positive majority when extracting an odd number $k$
of steps is 
\begin{equation}
P_{k}\left(y\right):=\sum_{h>k/2}\frac{k!}{h!\left(k-h\right)!}\,y^{h}\left(1-y\right)^{k-h}
\end{equation}
where the $h$ sum runs from $\left(k+1\right)/2$ to $k$. Then,
the urn function that describes an ERW with odd number $k>2$ of extractions
is 
\begin{equation}
\pi_{k}\left(y\right):=p\,P_{k}\left(y\right)+\left(1-p\right)\left(1-P_{k}\left(y\right)\right)=\left(1-p\right)+\left(2p-1\right)P_{k}\left(y\right)\label{eq:pik}
\end{equation}
it is a $k-$th degree polynomial, and is therefore non-linear for
all non-trivial values of the memory parameter $p$. 

\section{Strong convergence\protect\label{sec:Strong-convergence}}

\begin{figure}
\begin{centering}
\includegraphics[scale=0.29]{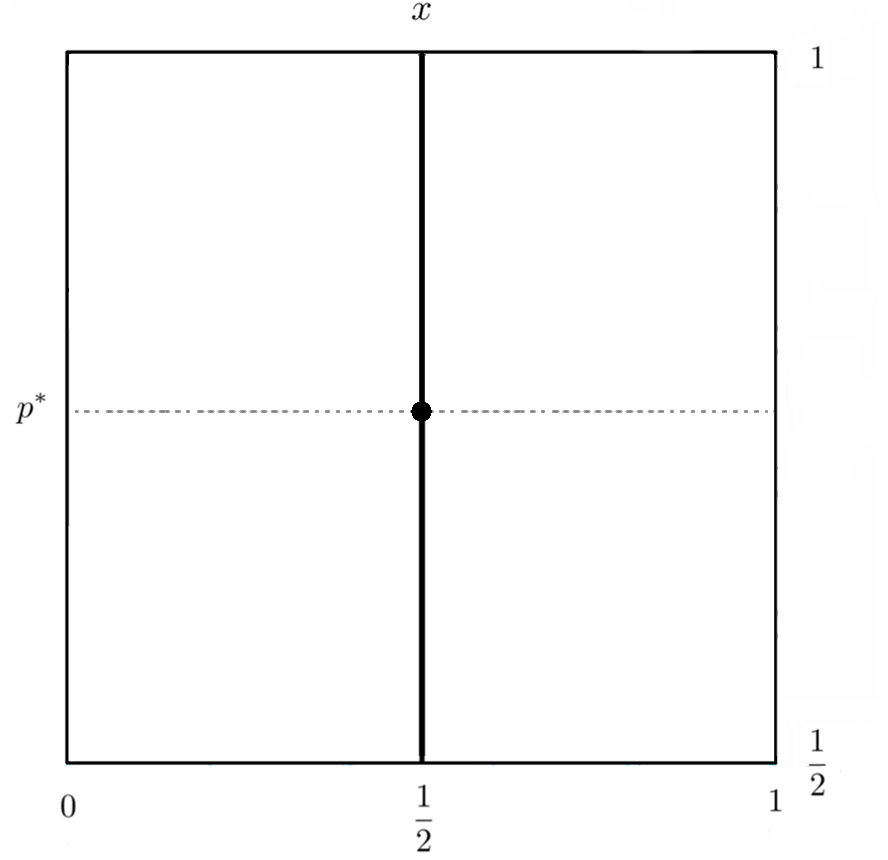}
\par\end{centering}
\caption{\protect\label{fig:3}Phase diagram $y$ vs $p$ for the entropy density
$\phi\left(y\right)$ of the classic ERW $k=1$, the diagram is shown
for $p>1/2$ only. For all $p<1$ the line at $y_{0}=1/2$ is the
only point of convergence for the density of black balls, although
there still is a critical $p^{*}$ above which the convergence of
$y_{N}$ is slowed according to the Pemantle's mechanism \cite{Jack=000020Harris,Pemantle=000020Touch,Franchini},
see Section \ref{sec:Large-deviations}. The line $y_{0}=1/2$ is
always a stable attractor for $y_{N}$ and the entropy is convex and
strictly negative in the whole region, except at the critical line
$y=1/2$, where it is zero. According to Eq. (\ref{eq:urnqe}), there
is a critical value at $p^{*}=3/4$ at which the derivative of the
urn function gets above $1/2$: for $p>p^{*}$ there is a shape change
in $\phi\left(y\right)$ in the neighborhood of $y=1/2$. \cite{FB}}
\end{figure}

The convergence properties of HLS urns for continuous urn functions
have been studied in \cite{HLS1,HLS2,Pemantle=000020Touch,Franchini,FB},
finding that the points of convergence of $y_{N}$ always belong to
the set of solutions of the equation (urn equation) 
\begin{equation}
\pi\left(y\right)=y,\label{eq:urnqe}
\end{equation}
and that these solutions are truly stable only if the derivative of
the urn function at those points is greater than one, ie if $\pi\left(y\right)$
crosses $y$ from top to bottom (down-crossing). For the classic ERW
(with $k=1$) the urn function $\pi_{1}\left(y\right)$ crosses $y$
at $1/2$ for any value of $p<1$, and therefore, $1/2$ is the only
possible point of convergence for the associated density of black
balls $y_{N}$, see Figure \ref{fig1}. This implies that the average
step $x_{N}$ converges to zero almost surely in the limit of $N\rightarrow\infty$
\begin{equation}
\lim_{N\rightarrow\infty}x_{N}=0,\ a.s.
\end{equation}
for all values of $p<1$ and with initial condition $x_{M}$, a phase
diagram for the ERW is shown in the Figure \ref{fig:3}. In the generalized
ERW with multiple extractions we see that a new phase appear below
some critical $p_{c}$. For $k=3$ already Eq. (\ref{eq:urnqe}) is
a third degree equation, and it can be solved with the well known
formula. In general, we find tree solutions: the first at $y_{0}=1/2$
and then \cite{FB}
\begin{equation}
y_{\pm}\left(p\right)=\frac{1}{2}\pm\sqrt{\frac{6p-5}{2p-1}},
\end{equation}
the term inside the square root is positive when $p\leq1/2$ and $p\geq5/6$,
but notice that $y_{\pm}\in\left[0,1\right]$ only if $p\in\left[1/2,1\right]$,
then, for $p$ below the critical value $p_{c}=5/6$ there is still
a unique stable solution at $1/2$ that crosses $y$ from top to bottom
(down crossing), see Figure \ref{fig:2}. \cite{FB}

\begin{figure}
\centering{}\includegraphics[scale=0.29]{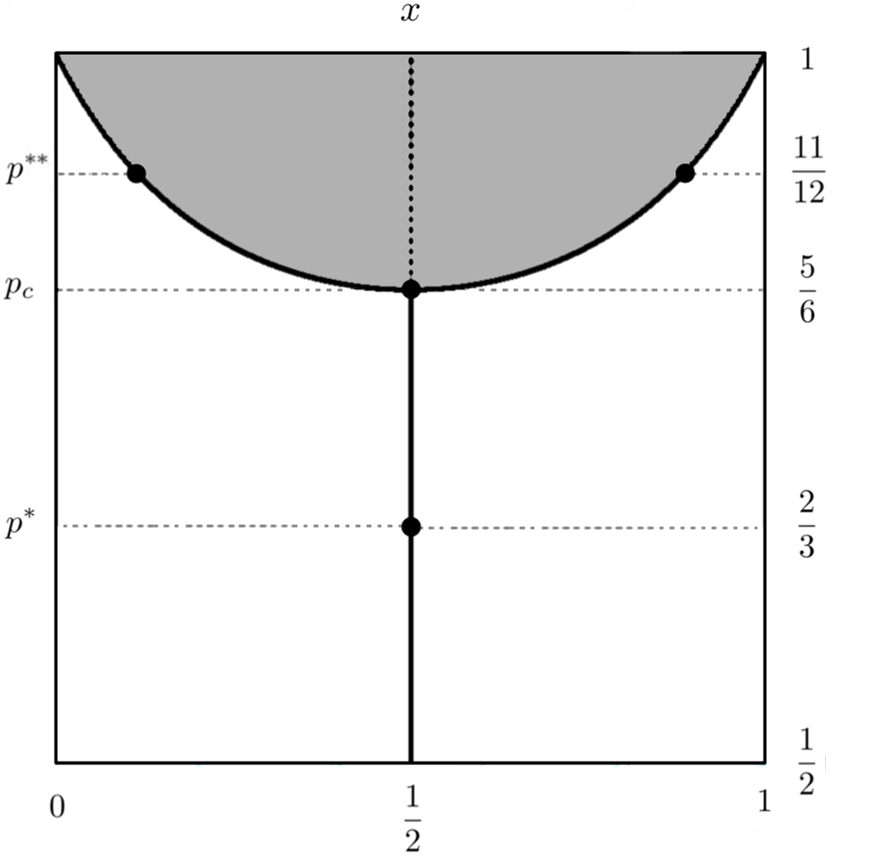}~\caption{\protect\label{fig:4}Phase diagram $x$ vs $p$ for the entropy density
$\phi\left(y\right)$ of the generalized ERW $k=3$. Above the critical
value $p_{c}=5/6$ the point $y_{0}=1/2$ becomes a unstable equilibrium,
and two new symmetric attractors arise according to Eq. (\ref{eq:ggoghuo}).
In the white colored region we still find a convex and negative $\phi\left(y\right)$,
except on the critical line, but notice that a new region appeared
above $p_{c}=5/6$, highlighted in darker shade, where $\phi\left(y\right)=0$,
ie the entropy is sub linear in $N$. The shape of $\phi\left(y\right)$
near the critical line is similar to the case $k=1$ for $p<p_{c}$,
except that here the point at which the derivative of the urn function
rise above $1/2$ is $p^{*}=2/3$. There is also another non-trivial
value $p^{**}=11/12$ at which the shape of $\phi\left(y\right)$
in the right (left) neighborhood of $y_{-}$($y_{+})$ has a change,
that is, when the derivative of the associated urn function in $y_{\pm}$
goes back below $1/2$. See also Figure \ref{fig:2}. \cite{FB}}
\end{figure}

Above $p_{c}$ the function still crosses $y$ at the point $1/2$,
but it now does from bottom to top, i.e. it is an up-crossing and
is therefore not stable. By the way, notice that for $p>p_{c}$ two
new solutions $y_{+}$ and $y_{-}$ also appear, that are both down-crossings
and can be stable attractors for $y_{N}$. Follows that for $p$ below
$p_{c}$ there are two attractors \cite{FB}
\begin{equation}
x_{\pm}\left(p\right)=\pm\sqrt{\frac{6p-5}{2p-1}},\label{eq:ggoghuo}
\end{equation}
separated by an unstable equilibrium $x_{0}=0$. The limit $x_{N}$
for $p>p_{c}$ is supported by 
\begin{equation}
\lim_{N\rightarrow\infty}x_{N}\in\left\{ x_{-}\left(p\right),\,x_{+}\left(p\right)\right\} ,\ a.s
\end{equation}
for any initial condition $x_{M}$. Since the urn functions that we
are considering never touch $0$ or $1$ for any $p\in\left(0,1\right)$,
for any initial condition $x_{M}$ that is fixed at finite time there
is a strictly positive probability to reach the nearby of any other
$x$ in a finite time at the very beginning of the process, then both
points $x_{\pm}$ carries some non-zero probability mass for any initial
condition at any finite time. Anyway it can be shown that the probability
mass near the point farther from $x_{M}$ will be exponentially suppressed
in $M$, and fixing the initial condition at some $M=o\left(N\right)$,
still divergent in $N$, will suppress one of the two possibilities,
and concentrate the probability mass in the attraction point $x_{\pm}$
that is closest to the initial $x_{M}$. The phases for $k=3$ are
shown in Figure \ref{fig:4}. 

Notice that in the limit of infinite $k$ the probability of finding
a majority of positive steps converges to the Heaviside theta function
\begin{equation}
P_{\infty}\left(y\right):=\theta\left(1-2y\right)
\end{equation}
and the urn function converges to a step function
\begin{equation}
\pi_{\infty}\left(y\right):=\left(1-p\right)+\left(2p-1\right)\theta\left(1-2y\right)
\end{equation}
that still crosses the diagonal at point $y_{0}=1/2$ (from top to
bottom) for $p<p_{c}=1/2$, and at $y_{-}=p$, $y_{+}=1-p$ if the
memory parameter is above $p_{c}$. Then also in the infinite $k$
limit there is a $p_{c}$ above which we find the same region of the
phase diagram that is observed for $k=3$. In fact, the phase diagram
shows the same structure for all $k>2$, apart from different $p_{c}$
and $x_{\pm}\left(p\right)$. For this reason, hereafter we restrict
our analysis to the cases $k=1$ and $k=3$ only. 

\section{Large deviations\protect\label{sec:Large-deviations}}

We can further refine the phase diagram by studying the entropy of
$y_{N}$ converging to $y$, 
\begin{equation}
\phi\left(y\right):=\lim_{N\rightarrow\infty}\frac{1}{N}\log P\left(y_{N}=\left\lfloor yN\right\rfloor /N\right),
\end{equation}
that is related to $\rho\left(x\right)$ by Eq. (\ref{eq:trans}):
\begin{equation}
\rho\left(x\right)=\phi\left(\left(1+x\right)/2\right).
\end{equation}
From the Large Deviations Theory of HLS urns (see the Corollary 5
of \cite{Franchini}), we know that, for any continuous and invertible
urn function, the limit $\phi\left(y\right)$ exists, that is strictly
convex, and is negative from $y=0$ up to the first point where the
urn function crosses the diagonal, is zero from there to the last
crossing, and then is convex negative again. 

Since for the classic ERW with $k=1$ the urn function is a line pivoting
around the center of the diagram, the CGF is zero only on the critical
point $y_{0}=1/2$ and strictly negative otherwise for $p<1$. For
the ERW with $k=3$: the picture below $p_{c}=5/6$ is qualitatively
the same, but above $p_{c}$ there is a region
\begin{equation}
y\in\left[y_{-}\left(p\right),\,y_{+}\left(p\right)\right]\label{eq:greyarea}
\end{equation}
where the entropy density is $\phi\left(y\right)=0$, i.e., the entropy
it sub-linear in $N$, see the dark colored region in Figure \ref{fig:4}.
In general, the main difference between the classic ERW and the other
cases ($k>2$) is the appearance of this novel phase above $p_{c}$,
where the entropy density is sub-linear in $N$ \cite{Franchini,FB}.

The asymptotic behavior of the entropy in this region could be inferred
by looking at the optimal trajectories taken by the process $y_{N}$
to reach some given limit value $y$. From the proof of Corollary
5 of \cite{Franchini}, we find that the entropy density of any HLS
urn with analytic and monotonic urn function $\pi$ is obtained trough
the following steps: let
\begin{equation}
\varphi:=\left\{ \varphi\left(\tau\right)\in\left[0,1\right]:\,0\leq\tau\leq1\right\} 
\end{equation}
and let $Q\left(y\right)$ be the set of absolutely continuous functions
on $\left[0,1\right]$ with initial value at zero, final value $y$,
and such that the derivative is between zero and one (Lipschitz),
\begin{equation}
Q\left(y\right):=\{\varphi\in C\left(\left[0,1\right]\right):\,0\leq\partial_{\tau}\varphi\left(\tau\right)\leq1,\,\varphi\left(0\right)=0,\,\varphi\left(1\right)=y\}.\label{eq:dfdfd}
\end{equation}
Also, define the auxiliary function
\begin{equation}
L\left(\alpha,\beta\right):=\alpha\log\left(\beta/\alpha\right)+\left(1-\alpha\right)\log\left(\left(1-\beta\right)/\left(1-\alpha\right)\right).\label{eq:sfsf}
\end{equation}
Then, in \cite{Franchini} it is shown that the entropy density $\phi\left(y\right)$
is obtained through solving the following variational problem in $\varphi$:
\begin{equation}
\phi\left(y\right)=-\inf_{\varphi\in Q\left(y\right)}I_{\pi}\left[\varphi\right]
\end{equation}
with rate function defined as follows: 
\begin{equation}
I_{\pi}\left[\varphi\right]:=\int_{0}^{1}d\tau\ L\left(\partial_{\tau}\varphi\left(\tau\right),\pi\left(\varphi\left(\tau\right)/\tau\right)\right),
\end{equation}
In Corollary 6 of \cite{Franchini} (see also \cite{FB,FranchiniHLS2025})
explicit formulas where given for the optimal trajectories of the
HLS model. From \cite{Franchini} we find that since $L$ is a negative
and concave function such that $L(\alpha,\beta)=0$ if and only if
$\alpha=\beta$, then when $I_{\pi}\left[\varphi\right]=0$ the trajectory
$\varphi$ must satisfy the following differential equation:
\begin{equation}
\partial_{\tau}\varphi\left(\tau\right)=\pi\left(\varphi\left(\tau\right)/\tau\right)
\end{equation}
with final condition $\varphi\left(1\right)=y$. Applying the substitution
\begin{equation}
\psi\left(\tau\right):=\varphi\left(\tau\right)/\tau,
\end{equation}
we find the differential equation of the scaling limit of the share:
\begin{equation}
\partial_{\tau}\psi\left(\tau\right)=\frac{\pi\left(\psi\left(\tau\right)\right)-\psi\left(\tau\right)}{\tau}
\end{equation}
with a final condition $\psi\left(1\right)=y$. The above equation
can be integrated, finding that the optimal strategy to achieve the
event $y_{N}=\left\lfloor yN\right\rfloor /N$ according to the Eq.
(\ref{eq:greyarea}) emanates from the closest unstable equilibrium
point. Moreover, from Corollary 6 of \cite{Franchini} follows the
uniqueness of the solutions, and that, given two trajectories $\psi_{1}\left(\tau\right)$
and $\psi_{2}\left(\tau\right)$, with $\psi_{1}\left(1\right)=y_{1}$,
$\psi_{2}\left(1\right)=y_{2}$, if $y_{1}<y_{2}$ then also $\psi_{1}\left(\tau\right)<\psi_{2}\left(\tau\right)$
for all $\tau\in\left(0,1\right)$.

Since any finite deviation from these trajectories has an exponential
cost on a time scale $O(N)$, the probability mass current can move
along these trajectories only, then the probability current passing
through $\left(\psi_{1}\left(\tau\right),\psi_{2}\left(\tau\right)\right)$
is a constant for all $\tau$, formally
\begin{equation}
P\left(y_{N}\in\left(y_{1},y_{2}\right)\right)=P\left(y_{\left\lfloor \tau N\right\rfloor }\in\left(\psi_{1}\left(\tau\right),\psi_{2}\left(\tau\right)\right)\right),
\end{equation}
and since the optimal trajectories can emanate only from the closest
unstable equilibrium, the probability of the event $y_{N}\in\left(y_{1},y_{2}\right)$
scales in $N$ like the entropy nearby that point. If $y_{0}$ is
that point, then a martingale analysis \cite{Franchini} would suggest
that 
\begin{equation}
P\left(y_{N}\in\left(y_{1},y_{2}\right)\right)={\textstyle O\,(N^{-\left(\chi-1\right)})}
\end{equation}
for large $N$, and with $\chi>1$, equal to the derivative of the
urn function in the point $y_{0}$. This analysis indicates that in
the sub-linear region the entropy is logarithmic in $N$. See also
in \cite{FB,Jack=000020LD} for a different estimation method.

\section{Cumulant Generating Function\protect\label{sec:Cumulant-Generating-Function}}

We can also study the behavior of $\phi\left(y\right)$ in the vicinity
of the critical line, by computing the Cumulant Generating Function
(CGF)
\begin{equation}
\zeta\left(\lambda\right):=\lim_{N\rightarrow\infty}\log\sum_{n\leq N}e^{-\lambda n}P\left(y_{N}=n/N\right),
\end{equation}
The right (left) behavior of the $\phi\left(y\right)$ near the critical
line can be deduced from the left (right) limit $\lambda\rightarrow0^{\pm}$
of the CGF before. Since the critical line is symmetrical around $y=1/2$
we only compute the limit from right, that is related to the entropy
density by the following Legendre Transform:
\begin{equation}
\phi\left(y\right)=\inf_{\lambda\in\left[0,\infty\right)}\left\{ \lambda y+\zeta\left(y\right)\right\} .
\end{equation}

In \cite{Franchini} is proven that if the urn function is monotonic
(invertible) the CGF satisfies the following nonlinear differential
equation: 
\begin{equation}
\partial_{\lambda}\,\zeta\left(\lambda\right)=\pi^{-1}\left({\textstyle \frac{e^{\,\zeta\left(\lambda\right)}-1}{e^{\,\lambda}-1}}\right)
\end{equation}
with $\pi^{-1}$ inverse urn function. Therefore, the linear urn function
for the classic ERW satisfies the differential equation 
\begin{equation}
\partial_{\lambda}\,\zeta\left(\lambda\right)=-\frac{1-p}{2p-1}+\frac{1}{2p-1}\left({\textstyle \frac{e^{\,\zeta\left(\lambda\right)}-1}{e^{\,\lambda}-1}}\right),
\end{equation}
the above equation can be integrated exactly: adapting the results
from Corollary 12 of \cite{Franchini} we find: \footnote{We where informed that Eq.s (2.44) and (2.45) Corollary 12 of \cite{Franchini}
have an inverted sign in front of the argument of the first exponential.
The same error propagated also to \cite{FB} and \cite{Franchini=000020Range=000020Urns}
(at least). I am grateful to Bernard Bercu, Michel Bonnefont and Adrien
Richou (Université Bordeaux) for spotting this error.}
\begin{multline}
1-e^{-\zeta\left(\lambda\right)}=\\
{\textstyle ={\textstyle \left(\frac{1-p}{2p-1}\right)}}e^{{\textstyle \left(\frac{1-p}{2p-1}\right)}\lambda}\left({\textstyle 1-e^{-\lambda}}\right)^{\left({\textstyle \frac{1}{2p-1}}\right)}\int_{1-e^{-\lambda}}^{1}dt\,\left(1-t\right)^{-\left({\textstyle {\textstyle \frac{3p-2}{2p-1}}}\right)}t^{{\textstyle -\left(\frac{1}{2p-1}\right)}}\label{eq:sssd}
\end{multline}
when $p>1/2$ and $\lambda>0$. Interestingly, in the region $p>1/2$
the CGF is never analytic at $\lambda=0$. Expanding the expression
before for small $\lambda$, we find a non vanishing term of $O(\lambda^{1/(2p-1)}\log\lambda)$
for $1/(2p-1)\in\mathbb{N}$, and $O(\lambda^{1/b})$ when $1/(2p-1)\in\mathbb{R}\setminus\mathbb{N}$,
i.e., derivatives $\left\lceil 1/(2p-1)\right\rceil $ and higher
are singular at zero \cite{Franchini}. 

Notice that, when $p>p^{*}=3/4$, i.e., when the derivative of the
urn function at the point of convergence $y_{0}$ climbs above $1/2$,
even the second order's cumulant diverges, then the shape $\phi\left(y\right)$
in the vicinity of $y_{0}=1/2$ is not Gaussian anymore for $p\in\left(p^{*},1\right)$,
see also Figure \ref{fig:3}. This indicates a phase change in the
convergence mechanism of $y_{N}$: below $p_{c}$, when the urn function
derivative at point $y_{0}$ is less than $1/2$, we expect for $y_{N}$
to cross the critical value infinitely many times in its evolution. 

But above $p_{c}$ the convergence of $y_{N}$ has a slowing down,
according to an interesting mechanism first described by Pemantle,
in \cite{Pemantle=000020Touch}, where $y_{N}$ approaches $y_{0}$
so slowly that it never crosses this point (almost surely), accumulating
in its right neighborhood. The convergence of the classic ERW in both
$p$ regimes has been further studied by Jack and Harris in \cite{Jack=000020Harris}.

Concerning the ERW with $k=3$ \cite{FB} the equation for the CGF
cannot be integrated with analytic methods (at best of our knowledge),
but looking at the behavior for small $\lambda$ we expect that below
the $p_{c}$ the same picture of the case with $k=1$ will arise,
with $p^{*}=5/6$, while above $p_{c}$ a new special value $p^{**}=11/12$
can be identified, where the derivative of the urn function at the
convergence point in $y_{+}\left(p\right)$ goes again below $1/2$.
We expect that in this last region the convergence mechanism of below
$p^{*}$ is restored, with different convergence point, see Figure
\ref{fig:4}. \cite{FB} 

\section{Relations with other models\protect\label{sec:Relations-with-other}}

Apart from the proposed generalization, the connection with the HLS
urn allows to put the ERW in relation with many models of physical
interest (and not only) that can be embedded into this model. 

For example, Brenig et al. \cite{Brenig} proved an isomorphism between
Quasi\textendash Polinomyal dynamical systems and multicolor balanced
urn processes. 

In reference \cite{Jack=000020LD}, Jack identified the urn function
that describes an interesting irreversible growth model introduced
by Klymko, Garrahan and Whitelam in \cite{KGW,KGGW}, the urn function
is written as follows:
\begin{equation}
\pi_{KGW}\left(y\right)=\frac{1+\tanh\left(J\left(2y-1\right)\right)}{2},
\end{equation}
controlled by the real parameter $J$. This urn function drives a
sub\textendash linear entropy region that is similar to the $k>2$
case of our generalized ERW model. 

Also, the HLS framework allows to connect the ERW with the classical
``Random Walk Range'' problem, \cite{Huges,Franchini=000020Range=000020Line,Franchini=000020Range=000020Urns,Franchini=000020Range}
that studies the number of different sites visited by a random walk
on the lattice $\mathbb{Z}^{d}$. In \cite{Franchini=000020Range=000020Urns}
we showed that the Range Problem can be exactly embedded into an HLS
model for some non\textendash linear urn function at any $d>1$. Although
for $d=2$ and $3$ a strongly non\textendash linear urn function
is observed, for $d\geq4$ the urn function gets surprisingly close
to a linear function in the Self-Avoiding Walk-like region of large
range values, that would be analogous to the classic ERW with $k=1$. 

Other than the physical models, one broad application of HLS urn is
as a simplified mechanism to explain the ``Lock\textendash In''
phenomenon in industrial and consumer behavior. \cite{Dosi=000020Ermoliev}
An influent model of Market Share between competing products that
can be modeled by the urn function $P_{k}\left(y\right)$ has been
introduced by Arthur in \cite{UM=000020Arthur}, and further developed
by other authors \cite{Dosi=000020Ermoliev,Franchini,FB,FranchiniHLS2025,Franchini=000020Range=000020Urns,Dosi=000020last,Gottfried,VanR,Gelast}.
In this model, two competing products gain customers according to
a majority mechanism, that correspond to the limit $p\rightarrow1$
of our generalized ERW with $k>2$. For an example of experimental
situations where this model actually comes into play we recommend
to give a look at the very interesting works by van de Rijt \cite{VanR}
and Gelastopulos et al. \cite{Gelast}, and also the Figure 8 of \cite{FranchiniHLS2025}
(provided by the same author).

\section*{Acknowledgements}

This paper serves as a proceeding for the homonymous talk in the session
CS19 Reinforcement Models: Elephant Random Walk, 44th Conference on
Stochastic Processes and their Applications, Wroclaw, July 14\textendash 18,
2025. It is a revision of an early version of \cite{FB}. I'm grateful
to Krishanu Maulik (Indian Statistical Institute) for inviting me.

\newpage{}

\end{document}